\documentclass[11pt]{article}
\usepackage{amsmath,amssymb,amsthm,amsxtra,graphicx,color,psfrag,epsfig,verbatim,mathrsfs,listings,fullpage}

\title{On pointless diagonal Fermat curves}
\author{Alexander P. McAvoy\footnote{Email:  \texttt{alexmcavoy@gmail.com}}}
\date{October 27, 2010}

\newtheorem{theorem}{Theorem}
\newtheorem{proposition}[theorem]{Proposition}

\newtheorem{lemma}[theorem]{Lemma}

\theoremstyle{definition}

\newtheorem{remark}{Remark}
\newtheorem{fact}{Fact}

\begin{document}

\maketitle

\begin{abstract}
We give an improved asymptotic upper bound on the number of diagonal Fermat curves $Ax^{\ell}+By^{\ell}=z^{\ell}$ over $\mathbb{F}_{q}$ with no $\mathbb{F}_{q}$-rational points, where $\ell$ is a prime number dividing $q-1$.
\end{abstract}

\setcounter{section}{-1}
\section{Introduction}

Among the most natural curves to consider over a finite field $\mathbb{F}_{q}$ are the diagonal Fermat curves, which are defined as the vanishing loci in $\mathbb{P}_{\mathbb{F}_{q}}^{2}$ of polynomials of the form $Ax^{\ell}+By^{\ell}-z^{\ell}$, where $A,B\in\mathbb{F}_{q}^{\times}$ and $\ell$ is a prime number.  However, such a curve, which we shall denote by $C_{q,\ell}\left(A,B\right)$, need not have any $\mathbb{F}_{q}$-rational points at all (for example, the curve $9x^{5}+4y^{5}=z^{5}$ has no $\mathbb{F}_{11}$-rational points).  Therefore, it makes sense to inquire as to when a given diagonal curve will in fact have any $\mathbb{F}_{q}$-rational points.  If $\ell$ does not divide $q-1$, then every element of $\mathbb{F}_{q}$ is an $\ell$th power, and we readily see that every diagonal curve has points over $\mathbb{F}_{q}$.  As a result, we shall assume that $\ell$ divides $q-1$ for the remainder of this paper.

To get some idea as to when the curve $C_{q,\ell}\left(A,B\right)$ has $\mathbb{F}_{q}$-rational points \textit{for every} $A,B\in\mathbb{F}_{q}^{\times}$, we first note that the requirement $A,B\neq 0$ implies that the projective variety $C_{q,\ell}\left(A,B\right)$ is smooth in $\mathbb{P}_{\mathbb{F}_{q}}^{2}$.  As $C_{q,\ell}\left(A,B\right)$ is a smooth algebraic curve of degree $\ell$, its genus is $g\left(C\right) =\left(\ell -1\right)\left(\ell -2\right) /2$.  Letting $N_{q,\ell}\left(A,B\right)$ denote the number of $\mathbb{F}_{q}$-rational points on $C_{q,\ell}\left(A,B\right)$ and letting
\begin{align*}
a_{q,\ell}\left(A,B\right) := N_{q,\ell}\left(A,B\right) -\left(q+1\right) ,
\end{align*}
we may apply the Hasse-Weil bound to see that $\left|a_{q,\ell}\left(A,B\right)\right|\leqslant 2g\left(C\right)\sqrt{q}=\left(\ell -1\right)\left(\ell -2\right)\sqrt{q}$.  We see immediately that if $N_{q,\ell}\left(A,B\right) =0$, then $q\ll\ell^{4}$.  Therefore, if $q\gg\ell^{4}$, then $N_{q,\ell}\left(A,B\right)\neq 0$ for every $A,B\in\mathbb{F}_{q}^{\times}$; in particular, every smooth diagonal curve has points.

However, in \cite{cohen}, Cohen observes experimentally that there are often values of $q\ll\ell^{4}$ such that $N_{q,\ell}\left(A,B\right)\neq 0$ for every $A,B\in\mathbb{F}_{q}^{\times}$, suggesting that the Hasse-Weil bound does not give an optimal asymptotic lower bound on $q$ guaranteeing points on every smooth diagonal curve.  Cohen then poses the following question:  If $Q\left(\ell\right)$ denotes the largest prime power $q$ such that $q\equiv 1\bmod\ell$ and $N_{q,\ell}\left(A,B\right) =0$ for some $A,B\in\mathbb{F}_{q}^{\times}$, can we improve on the trivial estimate of $Q\left(\ell\right)\ll\ell^{4}$?  Cohen's calculations that $Q\left(5\right) =11$, $Q\left(7\right) =71$, $Q\left(11\right) =419$, and $Q\left(13\right) =547$ suggest that the answer to his question is indeed in the affirmative.

The purpose of this paper is to study Cohen's question by looking at the \textit{average} value of $N_{q,\ell}\left(A,B\right)$ as we vary $A,B\in\mathbb{F}_{q}^{\times}$.  In order to state the main result of this paper, we first define
\begin{align*}
\mathscr{E}\left(q,\ell\right) := \left\{\left(A,B\right)\in\mathbb{F}_{q}^{\times}\times\mathbb{F}_{q}^{\times}\ :\ N_{q,\ell}\left(A,B\right) =0\right\} .
\end{align*}
Our main result is the following:

\begin{theorem}\label{moment2cor}
If $\delta$ is a real number between $0$ and $2$ and $q\gg\ell^{2+\delta}$, then $\#\mathscr{E}\left(q,\ell\right)\ll q^{2-\delta /\left(2+\delta\right)}$.
\end{theorem}

Since there are $\left(q-1\right)^{2}$ choices for $\left(A,B\right)\in\mathbb{F}_{q}^{\times}\times\mathbb{F}_{q}^{\times}$, we trivially have $\#\mathscr{E}\left(q,\ell\right)\ll q^{2}$, so Theorem $\ref{moment2cor}$ is indeed an improvement on the known result.  The upper bound of $\ell^{4}$ on $q$ in the hypothesis of Theorem $1$ comes from our previous comment about the Hasse-Weil bound since if $q\gg\ell^{4}$, then $\#\mathscr{E}\left(q,\ell\right) =0$.  If $0\leqslant\delta '\leqslant 3$ and $q\gg\ell^{1+\delta'}$, then our proof of Theorem $1$ yields
\begin{align*}
\#\mathscr{E}\left(q,\ell\right)\ll q^{3-\frac{2\delta '}{1+\delta '}} ,
\end{align*}
which clearly gives no improvement on the known bound $\#\mathscr{E}\left(q,\ell\right)\leqslant q^{2}$ for $0\leqslant\delta '\leqslant 1$, so we are only able to improve on the known result for $q$ in the range $\ell^{2}\ll q\ll\ell^{4}$.

Although Theorem $\ref{moment2cor}$ does not directly answer Cohen's question, our first proposition establishes $Q\left(\ell\right)\geqslant\ell^{2}$; in particular, only the range $\ell^{2}\ll q\ll\ell^{4}$ needs to be considered in order to provide a satisfactory answer to Cohen's question.

Our first lemma gives a closed formula for the sum $\sum_{A,B\in\mathbb{F}_{q}^{\times}}a_{q,\ell}\left(A,B\right)^{k}$ for $k=1,2$.  We find this closed formula by writing $N_{q,\ell}\left(A,B\right)$ explicitly in terms of Gauss and Jacobi sums.  Theorem $\ref{moment2cor}$ then follows as a corollary of the closed formula for $k=2$.  We subsequently find an explicit upper bound on
\begin{align*}
\left|\sum_{A,B\in\mathbb{F}_{q}^{\times}}a_{q,\ell}\left(A,B\right)^{k}\right|
\end{align*}
in terms of $q$ and $\ell$ for arbitrary $k\in\mathbb{N}$.

In the final section of this paper, we construct a multi-projective variety $\mathcal{V}\left(k\right)$ for every $k\in\mathbb{N}$ whose number of $\mathbb{F}_{q}$-rational points is related to the sum $\sum_{A,B\in\mathbb{F}_{q}^{\times}}a_{q,\ell}\left(A,B\right)^{k}$.  By proving that $\mathcal{V}\left(k\right)$ is rational for at least $k=1,2$, we establish once again the closed formula for the sum
\begin{align*}
\sum_{A,B\in\mathbb{F}_{q}^{\times}}a_{q,\ell}\left(A,B\right)^{2} .
\end{align*}

\section{Background, notation, and terminology}

By ``$k$th moment of the diagonal curve," we mean the value of the sum
\begin{align*}
\sum_{A,B\in\mathbb{F}_{q}^{\times}}a_{q,\ell}\left(A,B\right)^{k} .
\end{align*}
For a multiplicative character $\chi :\mathbb{F}_{q}^{\times}\rightarrow\mathbb{C}$ (which we shall extend to all of $\mathbb{F}_{q}$ by letting $\chi\left(0\right) =0$), the notation $g\left(\chi\right)$ denotes the Gauss sum $\sum_{\alpha\in\mathbb{F}_{q}}\chi\left(\alpha\right)\exp\left(2\pi i\alpha /p\right)$.  We denote by $\varepsilon$ the character with constant value $1$ on $\mathbb{F}_{q}^{\times}$.  We will briefly make reference to two Jacobi sums, defined as
\begin{align*}
J_{0}\left(\chi_{1},\dots ,\chi_{k}\right) &= \sum_{\alpha_{1}+\cdots +\alpha_{k}=0}\chi_{1}\left(\alpha_{1}\right)\cdots\chi_{k}\left(\alpha_{k}\right) ; \\
J\left(\chi_{1},\dots ,\chi_{k}\right) &= \sum_{\alpha_{1}+\cdots +\alpha_{k}=1}\chi_{1}\left(\alpha_{1}\right)\cdots\chi_{k}\left(\alpha_{k}\right) ,
\end{align*}
where $\alpha_{i}\in\mathbb{F}_{q}$.  We will need the following facts, whose proofs may be found in Chapter $10$ of \cite{rosen}:
\begin{fact}
$\left|g\left(\chi\right)\right| =\sqrt{q}$ if $\chi\neq\varepsilon$.
\end{fact}
\begin{fact}
If $\chi_{1},\dots ,\chi_{r}\neq\varepsilon$ but $\chi_{1}\cdots\chi_{r}=\varepsilon$, then $J_{0}\left(\chi_{1},\dots ,\chi_{r}\right) =\chi_{r}\left(-1\right)\left(q-1\right) J\left(\chi_{1},\dots ,\chi_{r-1}\right)$.
\end{fact}
\begin{fact}
If $\chi_{1},\dots ,\chi_{r},\chi_{1}\cdots\chi_{r}\neq\varepsilon$, then $g\left(\chi_{1}\right)\cdots g\left(\chi_{r}\right) =J\left(\chi_{1},\dots ,\chi_{r}\right) g\left(\chi_{1}\cdots\chi_{r}\right)$.
\end{fact}
\begin{fact}
The set of multiplicative characters of order $\ell$ on $\mathbb{F}_{q}$ forms a cyclic group of order $\ell$.
\end{fact}
\begin{fact}
If $\chi\neq\varepsilon$, then $\sum_{\alpha\in\mathbb{F}_{q}^{\times}}\chi\left(\alpha\right) =0$.  If $\chi =\varepsilon$, then $\sum_{\alpha\in\mathbb{F}_{q}^{\times}}\chi\left(\alpha\right) = q-1$.
\end{fact}
\begin{fact}
If $\chi\neq\varepsilon$, then $g\left(\chi\right) g\left(\overline{\chi}\right) =\chi\left(-1\right) q$.
\end{fact}

\section{Moments of the diagonal curve through Gauss sums}
In this section, we prove Theorem $\ref{moment2cor}$ by looking at the number of pointless diagonal Fermat curves \textit{on average}.  First, we prove a result establishing $Q\left(\ell\right)\gg\ell^{2}$.
\begin{proposition}
If $0<c<1$, then there exists a natural number $n\left(c\right)$ such that
\begin{align*}
n\left(c\right)\leqslant q\leqslant c\ell^{2} \implies \mathscr{E}\left(q,\ell\right)\neq\varnothing .
\end{align*}
\end{proposition}
\begin{proof}
We first note that the curve $Ax^{\ell}+By^{\ell}=z^{\ell}$ has no $\mathbb{F}_{q}$-rational points $\left[x:y:z\right]$ satisfying $xyz\neq 0$ if and only if $A$ is not a member of the zet $\left\{z^{\ell}-By^{\ell}\right\}_{y,z\in\mathbb{F}_{q}^{\times}}$ for some fixed $B\in\mathbb{F}_{q}^{\times}$ (we just take the chart $x=1$).  Since at most one of $x$, $y$, and $z$ may be equal to $0$ if $\left[x:y:z\right]$ is to be on a diagonal curve, we see that none of $A$, $B$, or $-A/B$ may be an $\ell$th power in $\mathbb{F}_{q}$ if the curve $Ax^{\ell}+By^{\ell}=z^{\ell}$ is to have no $\mathbb{F}_{q}$-rational points.  By fixing $B\in\mathbb{F}_{q}^{\times}$ such that $B$ is not an $\ell$th power, we see that the curve $Ax^{\ell}+By^{\ell}=z^{\ell}$ has no $\mathbb{F}_{q}$-rational points if and only if $A\not\in\left\{z^{\ell}-By^{\ell}\right\}_{y,z\in\mathbb{F}_{q}}$.

Now, if $0<c<1$ and $q\leqslant c\ell^{2}$, then the number of \textrm{pairs} of non-zero $\ell$th powers in $\mathbb{F}_{q}$ is $\left(q-1\right)^{2}/\ell^{2}\leqslant c\left(q-1\right)^{2}/q<cq$.  Therefore, if $q$ is sufficiently large (say $q\geqslant n\left(c\right)$ for some $n\left(c\right)\in\mathbb{N}$), the cardinality of the set $\left\{z^{\ell}-By^{\ell}\right\}_{y,z\in\mathbb{F}_{q}}$ is strictly less than $q-1$, which means that the desired value of $A$ can be found, completing the proof of the proposition.
\end{proof}

\begin{lemma}\label{moments12}
For the first and second moments of the diagonal curve, we have the closed formulae
\begin{align*}
\sum_{A,B\in\mathbb{F}_{q}^{\times}}a_{q,\ell}\left(A,B\right)^{k} = \begin{cases}0 & k=1, \\ q\left(q-1\right)^{2}\left(\ell -1\right)\left(\ell -2\right) & k=2.\end{cases}
\end{align*}
\end{lemma}
\begin{proof}
Letting $\chi_{1}$ and $\chi_{2}$ run over all non-trivial multiplicative characters of $\mathbb{F}_{q}$ of order $\ell$, we have
\begin{align*}
a_{q,\ell}\left(A,B\right) &= N_{q,\ell}\left(A,B\right) -\left(q+1\right) \\
&= \frac{1}{q-1}\sum_{\stackrel{\chi_{1},\chi_{2},\chi_{3}\neq\varepsilon}{\chi_{1}\chi_{2}\chi_{3}=\varepsilon}}\overline{\chi_{1}}\left(A\right)\overline{\chi_{2}}\left(B\right)\overline{\chi_{3}}\left(-1\right)J_{0}\left(\chi_{1},\chi_{2},\chi_{3}\right) \\
&= \sum_{\chi_{1},\chi_{2},\chi_{1}\chi_{2}\neq\varepsilon}\overline{\chi_{1}}\left(A\right)\overline{\chi_{2}}\left(B\right)J\left(\chi_{1},\chi_{2}\right) \\
&= \sum_{\chi_{1},\chi_{2},\chi_{1}\chi_{2}\neq\varepsilon}\overline{\chi_{1}}\left(A\right)\overline{\chi_{2}}\left(B\right)\left(\frac{g\left(\chi_{1}\right)g\left(\chi_{2}\right)}{g\left(\chi_{1}\chi_{2}\right)}\right)
\end{align*}
by Theorem $2$ on page $147$ of $\cite{rosen}$ and by Facts $1$ and $2$.  Therefore, the first moment of $a_{q,\ell}$ is
\begin{align*}
\sum_{A,B\in\mathbb{F}_{q}^{\times}}a_{q,\ell}\left(A,B\right) = \sum_{\chi_{1},\chi_{2},\chi_{1}\chi_{2}\neq\varepsilon}\left(\frac{g\left(\chi_{1}\right)g\left(\chi_{2}\right)}{g\left(\chi_{1}\chi_{2}\right)}\right)\sum_{A,B\in\mathbb{F}_{q}^{\times}}\overline{\chi_{1}}\left(A\right)\overline{\chi_{2}}\left(B\right) = 0
\end{align*}
by Fact $5$.  Similarly, the second moment of $a_{q,\ell}$ is
\begin{align*}
\sum_{A,B\in\mathbb{F}_{q}^{\times}}a_{q,\ell}\left(A,B\right)^{2} &= \sum_{\chi_{i1},\chi_{i2},\chi_{i1}\chi_{i2}\neq\varepsilon}\left(\frac{g\left(\chi_{11}\right)g\left(\chi_{12}\right)g\left(\chi_{21}\right)g\left(\chi_{22}\right)}{g\left(\chi_{11}\chi_{12}\right)g\left(\chi_{21}\chi_{22}\right)}\right) \\
&\quad\quad \times\sum_{A,B\in\mathbb{F}_{q}^{\times}}\overline{\chi_{11}\chi_{21}}\left(A\right)\overline{\chi_{12}\chi_{22}}\left(B\right) \\
&= \left(q-1\right)^{2}\sum_{\chi_{11},\chi_{12},\chi_{11}\chi_{12}\neq\varepsilon}\frac{g\left(\chi_{11}\right)g\left(\chi_{12}\right)g\left(\overline{\chi_{11}}\right)g\left(\overline{\chi_{12}}\right)}{g\left(\chi_{11}\chi_{12}\right)g\left(\overline{\chi_{11}\chi_{12}}\right)} \\
&= q\left(q-1\right)^{2}\left(\ell -1\right)\left(\ell -2\right)
\end{align*}
by Facts $5$ and $6$, which completes the proof.
\end{proof}

We are now in a position to prove Theorem $\ref{moment2cor}$:
\begin{proof}[Proof of Theorem \ref{moment2cor}]
By Lemma $\ref{moments12}$,
\begin{align*}
\sum_{A,B\in\mathbb{F}_{q}^{\times}}a_{q,\ell}\left(A,B\right)^{2} = q\left(q-1\right)^{2}\left(\ell -2\right)\left(\ell -2\right) \asymp q^{3}\ell^{2} .
\end{align*}
Therefore, for every $\left(A,B\right)\in\mathscr{E}\left(q,\ell\right)$, we have $a_{q,\ell}\left(A,B\right) =q+1$, so
\begin{align*}
q^{3}\ell^{2}\asymp \sum_{A,B\in\mathbb{F}_{q}^{\times}}a_{q,\ell}\left(A,B\right)^{2} \gg q^{2}\#\mathscr{E}\left(q,\ell\right) \implies \#\mathscr{E}\left(q,\ell\right)\ll q\ell^{2} .
\end{align*}
Letting $q\gg\ell^{2+\delta}$, we see that $\ell^{2}\ll q\ell^{-\delta}$ and $\ell\ll q^{1/\left(2+\delta\right)}$, hence
\begin{align*}
\#\mathscr{E}\left(q,\ell\right) \ll q\ell^{2} \ll q^{2}\ell^{-\delta} \ll q^{2-\delta /\left(2+\delta\right)} ,
\end{align*}
which completes the proof.
\end{proof}

\begin{figure}
\begin{center}
\caption{The graph of the the improved exponent given by Theorem $\ref{moment2cor}$ versus $\delta$ in the range $0<\delta <2$.  The improvement on the known result is clearly demonstrated by this figure.}
%
%
\begin{psfrags}%
\psfragscanon%
%
\psfrag{s03}[t][t]{\color[rgb]{0,0,0}\setlength{\tabcolsep}{0pt}\begin{tabular}{c}\textrm{The quantity }$\delta$\textrm{ coming from }$q\gg\ell^{2+\delta}$\end{tabular}}%
\psfrag{s04}[b][b]{\color[rgb]{0,0,0}\setlength{\tabcolsep}{0pt}\begin{tabular}{c}$2 - \delta / \left(2+\delta\right)$\textrm{, the improved exponent of Theorem }$\ref{moment2cor}$\end{tabular}}%
%
\psfrag{x01}[t][t]{0}%
\psfrag{x02}[t][t]{0.5}%
\psfrag{x03}[t][t]{1}%
\psfrag{x04}[t][t]{1.5}%
\psfrag{x05}[t][t]{2}%
%
\psfrag{v01}[r][r]{1.5}%
\psfrag{v02}[r][r]{1.55}%
\psfrag{v03}[r][r]{1.6}%
\psfrag{v04}[r][r]{1.65}%
\psfrag{v05}[r][r]{1.7}%
\psfrag{v06}[r][r]{1.75}%
\psfrag{v07}[r][r]{1.8}%
\psfrag{v08}[r][r]{1.85}%
\psfrag{v09}[r][r]{1.9}%
\psfrag{v10}[r][r]{1.95}%
\psfrag{v11}[r][r]{2}%
\resizebox{9cm}{!}{\includegraphics{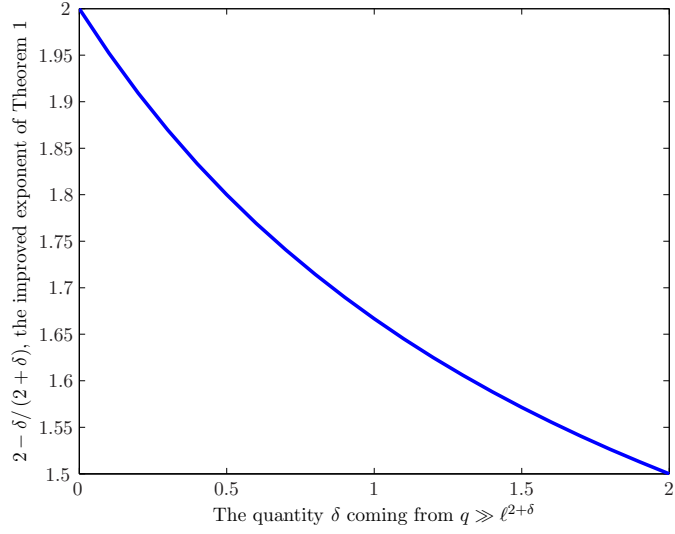}}
\end{psfrags}
\end{center}
\end{figure}

\begin{proposition}\label{momentbound}
For all $k\in\mathbb{N}$, we have
\begin{align*}
\left|\sum_{A,B\in\mathbb{F}_{q}^{\times}}a_{q,\ell}\left(A,B\right)^{k}\right| \leqslant q^{k/2}\left(q-1\right)^{2}\left(\ell -1\right)\left(\ell -2\right)^{k-1}\left(\frac{\left(\ell -1\right)^{k-1}-\left(-1\right)^{k-1}}{\ell}\right) .
\end{align*}
\end{proposition}
\begin{proof}
Lemma $\ref{moments12}$ proves this result for $k=1,2$.  Let $k\geqslant 3$ be an integer.  Then,
\begin{align*}
\sum_{A,B\in\mathbb{F}_{q}^{\times}}a_{q,\ell}\left(A,B\right)^{k} &= \sum_{A,B\in\mathbb{F}_{q}^{\times}}\left(\sum_{\chi_{1},\chi_{2},\chi_{1}\chi_{2}\neq\varepsilon}\overline{\chi_{1}}\left(A\right)\overline{\chi_{2}}\left(B\right)\left(\frac{g\left(\chi_{1}\right)g\left(\chi_{2}\right)}{g\left(\chi_{1}\chi_{2}\right)}\right)\right)^{k} \\
&= \left(q-1\right)^{2}\sum_{i=1}^{k-1}\sum_{\chi_{i,1},\chi_{i,2},\chi_{i,1}\chi_{i,2}\neq\varepsilon}\prod_{j=1}^{k-1}\left(\frac{g\left(\chi_{j,1}\right)g\left(\chi_{j,2}\right)}{g\left(\chi_{j,1}\chi_{j,2}\right)}\right) \\
&\quad\quad \times\left(\frac{g\left(\overline{\chi_{1,1}\cdots\chi_{k-1,1}}\right)g\left(\overline{\chi_{1,2}\cdots\chi_{k-1,2}}\right)}{g\left(\overline{\left(\chi_{1,1}\chi_{1,2}\right)\cdots\left(\chi_{k-1,1}\chi_{k-1,2}\right)}\right)}\right) \\
&= \frac{\left(q-1\right)^{2}}{q^{k-1}}\sum_{\psi_{1},\dots ,\psi_{k-1},\psi_{1}\cdots\psi_{k-1}\neq\varepsilon}\psi_{1}\cdots\psi_{k-1}\left(-1\right)\left(\frac{g\left(\overline{\psi_{1}}\right)\cdots g\left(\overline{\psi_{k-1}}\right)}{g\left(\overline{\psi_{1}\cdots\psi_{k-1}}\right)}\right) \\
&\quad\quad \times\sum_{i=1}^{k-1}\sum_{\chi_{i,1}\chi_{i,2}=\psi_{i}}\Big(g\left(\chi_{1,1}\right)g\left(\chi_{1,2}\right)\Big)\cdots\Big(g\left(\chi_{k-1,1}\right)g\left(\chi_{k-1,2}\right)\Big) \\
&\quad\quad \times g\left(\overline{\chi_{1,1}\cdots\chi_{k-1,1}}\right)g\left(\overline{\chi_{1,2}\cdots\chi_{k-1,2}}\right) .
\end{align*}
Now, each sum $\sum_{\chi_{i,1}\chi_{i,2}=\psi_{i}}$ is clearly of size $\left(\ell -2\right)^{k-1}$.  Letting $S\left(k-1\right)$ denote
\begin{align*}
\#\left\{\psi_{1},\dots ,\psi_{k-1},\psi_{1}\cdots\psi_{k-1}\neq\varepsilon\right\} ,
\end{align*}
we see that $S\left(k-1\right) = \left(\ell-1\right) S\left(k-3\right) +\left(\ell -2\right) S\left(k-2\right)$, so $S\left(k-1\right) =\left(\ell -1\right)^{k-1}-S\left(k-2\right)$.  Therefore, inductively, we obtain the formula
\begin{align*}
S\left(k-1\right) = \left(-1\right)^{k-1}\sum_{i=1}^{k-1}\left(1-\ell\right)^{i} = \left(-1\right)^{k-1}\left(\ell -1\right)\left(\frac{\left(1-\ell\right)^{k-1}-1}{\ell}\right)
\end{align*}
Recalling that $\left|g\left(\chi\right)\right| = \sqrt{q}$ for $\chi\neq\varepsilon$  (Fact $1$) and applying the triangle inequality, we find that
\begin{align*}
\left|\sum_{A,B\in\mathbb{F}_{q}^{\times}}a_{q,\ell}\left(A,B\right)^{k}\right| \leqslant q^{k/2}\left(q-1\right)^{2}\left(\ell -2\right)^{k-1}\left(-1\right)^{k-1}\left(\ell -1\right)\left(\frac{\left(1-\ell\right)^{k-1}-1}{\ell}\right) ,
\end{align*}
which gives the desired result.
\end{proof}

\begin{table}\label{mytable}
\begin{center}
\caption{Numerical data illustrating Proposition $\ref{momentbound}$ up to the 10th moment of $a_{q}$ using $805$ pairs $\left(p,\ell\right)$ in the range $211\leqslant p\leqslant 5233$ and $13\leqslant\ell\leqslant 2459$.}
\begin{tabular}{|c|c|c|}
\hline
\textbf{k} & \textbf{Bound given by Proposition $\ref{momentbound}$} & \textbf{Experimental bound using $805$ pairs $\left(p,\ell\right)$} \\
\hline
$2$ & $p^{3}\ell^{2}$ & $p^{3.0182}\ell^{2.6808}$ \\
\hline
$3$ & $p^{3.5}\ell^{4}$ & $p^{3.4580}\ell^{3.9892}$ \\
\hline
$4$ & $p^{4}\ell^{6}$ & $p^{4.0510}\ell^{5.8134}$ \\
\hline
$5$ & $p^{4.5}\ell^{8}$ & $p^{4.5160}\ell^{7.3777}$ \\
\hline
$6$ & $p^{5}\ell^{10}$ & $p^{5.0893}\ell^{9.1752}$ \\
\hline
$7$ & $p^{5.5}\ell^{12}$ & $p^{5.5723}\ell^{10.8686}$ \\
\hline
$8$ & $p^{6}\ell^{14}$ & $p^{6.1296}\ell^{12.6346}$ \\
\hline
$9$ & $p^{6.5}\ell^{16}$ & $p^{6.6169}\ell^{14.3598}$ \\
\hline
$10$ & $p^{7}\ell^{18}$ & $p^{7.1711}\ell^{16.1363}$ \\
\hline
\end{tabular}
\end{center}
\end{table}

It is natural to now ask if we can improve upon Theorem $\ref{moment2cor}$ using higher moments of $a_{q}$.  Using the bound given by Proposition $\ref{momentbound}$ and the method used to prove Lemma $\ref{moments12}$, we see that
\begin{align*}
q^{\left(k+4\right) /2}\ell^{2k-2} \gg \sum_{A,B\in\mathbb{F}_{q}^{\times}}a_{q,\ell}\left(A,B\right)^{k} \gg \sum_{\left(A,B\right)\in\mathscr{E}\left(q,\ell\right)}q^{k} = q^{k}\#\mathscr{E}\left(q,\ell\right) .
\end{align*}
Therefore, we have $\#\mathscr{E}\left(q,\ell\right)\ll q^{\left(4-k\right) /2}\ell^{2k-2}$.  Setting $q\gg\ell^{2+\delta}$ for $0<\delta <2$, we obtain
\begin{align*}
\#\mathscr{E}\left(q,\ell\right)\ll q^{\left(k+4\right) /2}\ell^{-2-k\delta}\ll q^{\left(k+4\right) /2-\left(2+k\delta\right) /\left(2+\delta\right)} .
\end{align*}
It is easily checked that for any $\delta$ in the range $0<\delta <2$,
\begin{align*}
\frac{k+4}{2}-\frac{2+k\delta}{2+\delta}\leqslant 2-\frac{\delta}{2+\delta} \implies k\leqslant 2 ,
\end{align*}
so we would need to improve upon the bound given in Proposition $\ref{momentbound}$ to improve Theorem $\ref{moment2cor}$.

\begin{remark}
Notice that the formulae for the first and second moments are polynomials in $q$ and $\ell$.  In general, we cannot hope for such a nice closed formula for higher moments.  However, given the Gauss sum expansion used in the proof of Theorem $\ref{moment2cor}$, we see that it is plausible that such a closed formula, should it exist, might be a polynomial in $q^{1/2}$ and $\ell$ (the identity of such a closed formula will be explored more in the final section).  The following theorem shows that, in fact, no polynomial in $q^{1/2}$ and $\ell$ can be equal to the third or fourth moment of $a_{q}$:
\end{remark}

\begin{theorem}\label{noclosed}
There does not exist a polynomial $f\left(x,y\right)\in\mathbb{Z}\left[x^{1/2},y\right]$ such that
\begin{align*}
f\left(q,\ell\right) = \sum_{A,B\in\mathbb{F}_{q}^{\times}}a_{q,\ell}\left(A,B\right)^{k}
\end{align*}
for every $q\in\mathbb{F}_{q}$ and $\ell$ prime satisfying $q\equiv 1\bmod\ell$ if $k=3,4$.
\end{theorem}
\begin{proof}
Given the method of proof used in Proposition $\ref{momentbound}$, we see that should such an $f$ exist, it would be of degree at most $k+4$ in $x^{1/2}$ and of degree at most $2k-2$ in $y$.  Therefore, suppose that such a polynomial exists and write
\begin{align*}
f\left(x,y\right) = \sum_{\stackrel{0\leqslant i\leqslant k+4}{0\leqslant j\leqslant 2k-2}}c_{i,j}x^{i/2}y^{j} .
\end{align*}
Should such a polynomial be equal to the $k$th moment of $a_{q}$ for an arbitrary prime power $q$, it will certainly be true for $q=p$.  Therefore, letting $z_{k}\left(p,\ell\right)$ denote $\sum_{A,B\in\mathbb{F}_{q}^{\times}}a_{q,\ell}\left(A,B\right)^{k}$, we can choose $\left(k+5\right)\left(2k-1\right)$ pairs $\left(p,\ell\right)$ and a square matrix $M$ such that
\begin{align*}
\begin{pmatrix}z_{k}\left(p_{1},\ell_{1}\right) \\ \vdots \\ z_{k}\left(p_{\left(k+5\right)\left(2k-1\right)}\ell_{\left(k+5\right)\left(2k-1\right)}\right)\end{pmatrix} = M\begin{pmatrix}c_{0,0} \\ c_{1,0} \\ \vdots \\ c_{k+4,0} \\ c_{0,1} \\ c_{1,1} \\ \vdots \\ c_{k+4,2k-2}\end{pmatrix} .
\end{align*}
For $k=3$, we choose two distinct sets $S_{1}$ and $S_{2}$ each containing $40$ pairs of $\left(p,\ell\right)$ such that the resulting matrices $M_{1}$ and $M_{2}$ (corresponding to $M$ above) are invertible.  Using these sets, we calculate the corresponding third moments $\left\{z_{3}\left(p,\ell\right)\right\}_{\left(p,\ell\right)\in S_{i}}$ for $i=1,2$ (call the resulting third moment vectors $v_{1}$ and $v_{2}$, respectively).  Using PARI/GP, we calculate that $M_{1}^{-1}v_{1}\neq M_{2}^{-1}v_{2}$.  Clearly, if the desired polynomial $f$ should exist, we would have to have $M_{1}^{-1}v_{1}=M_{2}^{-1}v_{2}$, so such a closed formula for the third moment indeed does not exist.  The same calculation (but with two sets of $54$ pairs of $\left(p,\ell\right)$ instead of $40$ such pairs) gives the result for the fourth moment.
\end{proof}
The closed formula for the second moment is easily verified numerically using the method outlined in the above proof.  Although the above proof provides a definitive answer to the question of the existence of a polynomial in $q^{1/2}$ and $\ell$ that is equal to the third or fourth moment, it is not computationally feasible to use in order to answer this question for higher moments of $a_{q,\ell}$.

It is natural to ask whether or not we can find closed formulae for higher moments of the diagonal curve.  By examining the proof of Lemma $\ref{moments12}$, we see that the fact that a closed formula for the second moment can be found using Gauss sums seems fortuitous:  Each Gauss sum could be paired with its conjugate and subsequently canceled.  However, for the third and higher moments, we do not find such cancellation of Gauss sums.  We therefore move to the next section in which we introduce a geometric interpretation of the higher moments.

\section{Moments of the diagonal curve through geometry}

In this section, we let $C_{q,\ell}\left(A,B,C\right)$ denote the vanishing locus of the equation $Ax^{\ell}+By^{\ell}+Cz^{\ell}$ in $\mathbb{P}_{\mathbb{F}_{q}}^{2}$ for $\left[A:B:C\right]\in\mathbb{P}_{\mathbb{F}_{q}}^{2}$ and we let $N_{q,\ell}\left(A,B,C\right)$ be the number of $\mathbb{F}_{q}$-rational points on $C_{q,\ell}\left(A,B,C\right)$.  We can actually look at an entire family of diagonal curves by defining $\mathcal{V}\left(1\right)$ to be the vanishing locus of $Ax^{\ell}+By^{\ell}+Cz^{\ell}$ in $\mathbb{P}_{\mathbb{F}_{q}}^{2}\times\mathbb{P}_{\mathbb{F}_{q}}^{2}$.  Explicitly, we have
\begin{align*}
\mathcal{V}\left(1\right) := \left\{\left(\left[A:B:C\right] ,\left[x:y:z\right]\right)\in\mathbb{P}_{\mathbb{F}_{q}}^{2}\times\mathbb{P}_{\mathbb{F}_{q}}^{2}\ :\ Ax^{\ell}+By^{\ell}+Cz^{\ell}=0\right\} .
\end{align*}
The key observation of this section is that the value of the sum
\begin{align*}
\sum_{\left[A:B:C\right]\in\mathbb{P}_{\mathbb{F}_{q}}^{2}}N_{q,\ell}\left(A,B,C\right)^{k}
\end{align*}
is equal to the number of $\mathbb{F}_{q}$-rational points on the $k$-fold fibre product
\begin{align*}
\mathcal{V}\left(k\right) := \mathcal{V}\left(1\right) \times_{\mathbb{P}_{\mathbb{F}_{q}}^{2}} \cdots \times_{\mathbb{P}_{\mathbb{F}_{q}}^{2}} \mathcal{V}\left(1\right) ,
\end{align*}
where the fibre product is being taken over the $\left[A:B:C\right]$-coefficient space.

It is more intuitive and natural to arrive at Theorem $\ref{moment2cor}$ by looking at the geometry of the moments of diagonal curve rather than explicit Gauss sum calculations.  Since $a_{q,\ell}\left(A,B\right)$ and $N_{q,\ell}\left(A,B\right)$ differ by $\left(q+1\right)$, it follows that we can find the $k$th moment of $a_{q,\ell}\left(A,B\right)$ if and only if we can determine the $k$th moment of $N_{q,\ell}\left(A,B\right)$.  In the previous section, it was more convenient to look at moments of $a_{q,\ell}\left(A,B\right)$ due to their nice expansion in terms of Gauss sums.  Denoting by $\mathcal{V}_{0}\left(k\right)$ the smooth part of the multi-projective variety $\mathcal{V}\left(k\right)$, which is precisely the part of the fibre product taken over the set of points $\left[A:B:C\right]\in\mathbb{P}_{\mathbb{F}_{q}}^{2}$ such that $ABC\neq 0$, we have:

\begin{proof}[Geometric proof of Lemma \ref{moments12}]
We just demonstrate that $\mathcal{V}\left(1\right)$ and $\mathcal{V}\left(2\right)$ are rational and subsequently count the number of points on the fibre product.  Since we are looking to prove Lemma $\ref{moments12}$, which only deals with $\mathcal{V}_{0}\left(1\right)$ and $\mathcal{V}_{0}\left(2\right)$, we may assume that $C=-1$ and $A,B\neq 0$:

\noindent\textbf{First moment:}  The derivation of the first moment is nearly trivial:  If $x=0$, there are $\left(q-1\right)^{2}$ points on the curve.  If $x\neq 0$, we can solve for $A$ to see that there are $q\left(q-1\right)^{2}$ points on the curve, for a total of $\left(q-1\right)^{2}\left(q+1\right)$ points, from which it follows immediately that
\begin{align*}
\sum_{A,B\in\mathbb{F}_{q}^{\times}}a_{q,\ell}\left(A,B\right) = \# \mathcal{V}_{0}\left(1\right) -\left(q-1\right)^{2}\left(q+1\right) = \left(q-1\right)^{2}\left(q+1\right) - \left(q-1\right)^{2}\left(q+1\right) = 0.
\end{align*}

\noindent\textbf{Second moment:}  We simply set up the system
\begin{align*}
\begin{pmatrix}x_{1}^{\ell} & y_{1}^{\ell} \\ x_{2}^{\ell} & y_{2}^{\ell}\end{pmatrix}\begin{pmatrix}A \\ B\end{pmatrix} = \begin{pmatrix}z_{1}^{\ell} \\ z_{2}^{\ell}\end{pmatrix}
\end{align*}
and consider the cases $z_{1}=z_{2}=0$; $z_{1}=1$ and $z_{2}=0$; and $z_{1}=z_{2}=1$, which are the only cases to consider by symmetry.  Simple manipulation using linear algebra shows that there are $\left(q-1\right)^{2}\ell$ points in the first case, $\left(q-1\right)^{2}\left(q-\ell +1\right)$ points in the second case, and
\begin{align*}
\left(q-1\right)^{2}\left(\left(q-1\right)^{2}+4q-\left(q-1\right)\ell -2-2\ell\left(q-\ell +1\right)\right)
\end{align*}
points in the third case.  Therefore, we find that the second moment is
\begin{align*}
\sum_{A,B\in\mathbb{F}_{q}^{\times}}a_{q,\ell}\left(A,B\right)^{2} &= \#\mathcal{V}_{0}\left(2\right) - 2\left(q+1\right) \#\mathcal{V}_{0}\left(1\right) + \left(q-1\right)^{2}\left(q+1\right)^{2} \\
&= \left(q-1\right)^{2}\left(q+1\right)^{2} + q\left(q-1\right)^{2}\left(\ell -1\right)\left(\ell -2\right) - \left(q-1\right)^{2}\left(q+1\right)^{2} \\
&= q\left(q-1\right)^{2}\left(\ell -1\right)\left(\ell -2\right) ,
\end{align*}
which completes the proof of Lemma $\ref{moments12}$.
\end{proof}

\section{Future work}

The singular locus of $\mathcal{V}\left(k\right)$ is contained in the portion of the variety above the triangle $ABC=0$ in $\mathbb{P}_{\mathbb{F}_{q}}^{2}$.  By resolving these singularities, we would have a smooth variety $\widetilde{\mathcal{V}}\left(k\right)$ birational to $\mathcal{V}\left(k\right)$, and we could subsequently express the number of points on $\widetilde{\mathcal{V}}\left(k\right)$ in terms of the Lefschetz fixed point formula by the Weil conjectures.  By determining $\#\widetilde{\mathcal{V}}\left(k\right) -\#\mathcal{V}\left(k\right)$, we may be able to use Proposition $\ref{momentbound}$ to show that some of the cohomology groups of $\widetilde{\mathcal{V}}\left(k\right)$ are trivial.  Indeed, the dimension of $\widetilde{\mathcal{V}}\left(k\right)$ must be $k+2$ since that is the dimension of $\mathcal{V}\left(k\right)$, and for $\lambda$ coprime to $q$, we have
\begin{align*}
\left|\#\widetilde{\mathcal{V}}\left(k\right)\right| &= \left|\sum_{i=0}^{2\left(k+2\right)}\left(-1\right)^{i}\textrm{tr}\left(\textrm{Frob}_{q}\ |\ H^{i}\left(\overline{\widetilde{\mathcal{V}}\left(k\right)},\mathbb{Q}_{\lambda}\right)\right)\right| \\
&\leqslant \sum_{i=0}^{2\left(k+2\right)}\left(\dim H^{i}\left(\overline{\widetilde{\mathcal{V}}\left(k\right)},\mathbb{Q}_{\lambda}\right)\right)q^{i/2}
\end{align*}
since the magnitude of each eigenvalue of the action of the Frobenius on the $i$th cohomology group is $q^{i/2}$ (see Appendix C of \cite{hart}).  Therefore, if the varieties $\mathcal{V}\left(k\right)$ are properly desingularized and the pullback of the singular locus is calculated, the bound of $\left(k+4\right) /2$ on the exponent of $q$ in Proposition $\ref{momentbound}$ might give some insight into the identity of the first few cohomology groups since the Lefschetz fixed point theorem yields a bound of $k+2$ on the exponent of $q$.

In practice, actually desingularizing these varieties has proven to be quite difficult, not to mention actually calculating the pullback of the singular locus under desingularization.

\section*{Acknowledgments}

I am greatly indebted to Siman Wong for proposing this problem and for all of his suggestions during the course of the research that resulted in this paper (which serves as my undergraduate thesis).  I would also like to thank Hans Johnston for providing the computers on which the PARI/GP calculations were carried out, as well as Jenia Tevelev, Paul Hacking, and Jessica Sidman for helpful conversations about the geometry of the diagonal varieties.  Finally, I would like to thank Farshid Hajir for his comments and suggestions on earlier drafts of this paper.

\bibliographystyle{amsalpha}

\end{document}